\documentclass{amsart}

\usepackage{amssymb}
\usepackage{amsmath}
\usepackage[all]{xy}
\usepackage{enumerate}
\usepackage{hyperref}

\newtheorem{theorem}[subsection]{Theorem}
\newtheorem{cor}[subsection]{Corollary}
\newtheorem{lemma}[subsection]{Lemma}

\newtheorem{conj}[subsection]{Conjecture}

\theoremstyle{remark}
\newtheorem{exa}[subsection]{Example}
\newtheorem{rem}[subsection]{Remark}

\theoremstyle{definition}

\newcommand{\Z}{{\mathbb Z}}
\newcommand{\R}{{\mathbb R}}
\DeclareMathOperator{\E}{E}
\DeclareMathOperator{\U}{U}

\begin{document}

\title{On the Log-Concavity of Hilbert Series of Veronese Subrings and Ehrhart Series}

\author{Matthias Beck}
\address{Department of Mathematics\\
         San Francisco State University\\
         San Francisco, CA 94132\\
         U.S.A.}
\email{beck@math.sfsu.edu}

\author{Alan Stapledon}
\address{Department of Mathematics, University of Michigan, Ann Arbor, MI 48109, U.S.A.}
\email{astapldn@umich.edu}

\date{13 October 2008}

\thanks{The authors would like to thank Alexander Barvinok, Jesus De Loera, Sergey Fomin, Joseph Gubeladze, Mircea Musta\c t\v a, Sam Payne, John Stembridge, Volkmar Welker, and an anonymous referee for useful discussions and helpful suggestions.
% They thank Volkmar Welker for bringing \cite{BWVeronese} to their attention.
The first author was partially supported by the NSF (research grant DMS-0810105), and the second author was partially supported by an Eleanor Sophia Wood
travelling scholarship from the University of Sydney.}

% 52B20 Lattice polytopes
% 05A15 Exact enumeration problems, generating functions
% 13C14 Cohen-Macaulay module

\maketitle

\begin{abstract}

For every positive integer $n$, consider the linear
operator $\U_{n}$ on polynomials of degree at most $d$
with integer coefficients defined as follows: if we
write $\frac{h(t)}{(1 - t)^{d + 1}} = \sum_{m \geq 0}
g(m) \, t^{m}$, for some polynomial $g(m)$ with rational
coefficients, then  $\frac{\U_{n}h(t)}{(1- t)^{d + 1}}
= \sum_{m \geq 0} g(nm) \, t^{m}$. We show that there
exists a positive integer $n_{d}$, depending only on
$d$, such that if $h(t)$ is a polynomial  of degree at
most $d$ with nonnegative integer coefficients and
$h(0) = 1$, then for $n \geq n_{d}$, $\U_{n}h(t)$
has simple, real, negative roots and positive,
strictly log concave and strictly unimodal
coefficients. Applications are given to Ehrhart
$\delta$-polynomials and unimodular triangulations of
dilations of lattice polytopes, as well as Hilbert
series of Veronese subrings of Cohen--Macauley graded
rings.

%define $\U_{n}h(t)$ to be the polynomial of degree at most $d$ with integer coefficients satisfying
%$\frac{\U_{n}h(t)}{(1- t)^{d + 1}} =
%\sum_{m \geq 0} g(nm) \, t^{m}$.

%Let $N$ be a lattice and $P \subset N \otimes_{\mathbb{Z}} \mathbb{R}$ a lattice polytope, i.e., the convex hull of finitely many points in $N$. Ehrhart's theorem asserts that the lattice-point counting function $f_{P}(m) := \# \left( mP \cap N \right)$ is a polynomial, and thus $1 + \sum_{m \geq 1} f_{P}(m) \, t^{m} = \frac{ \delta_{P}(t) }{ (1 - t)^{d + 1} }$ for some polynomial $\delta_{P}(t)$, the \emph{$\delta$-vector} of $P$. Motivated by the Knudsen--Mumford--Waterman Conjecture about the existence of unimodular triangulations for all sufficiently large dilates of $P$ and a recent theorem of Athanasiadis--Hibi--Stanley asserting that such triangulations imply certain inequalities for the $\delta$-vector, we study $\delta_{ nP } (t) = \delta_{0}(n) + \delta_{1}(n) \, t + \cdots +  \delta_{d}(n) \, t^{d}$ as $n$ grows. We prove that for sufficiently large $n$ (with a bound depending only on $\dim P$), $\delta_{0}(n) < \delta_{d}(n) < \delta_{1}(n) < \cdots < \delta_{i}! (n) < \delta_{d - i}(n) < \delta_{i + 1}(n) < \cdots < \delta_{\lfloor \frac{d + 1}{2} \rfloor}(n)$; in particular, $\delta_{ nP } (t)$ is unimodal.
\end{abstract}

%%%%%%%%%%%%%%%%%%%%%%%%%%%%%%%%%%%%%%%%%%%%%%%%%%%%%%%%%%%%%%%%%%%%%%%%%%%%%%%%

\section{Introduction}

Fix a positive integer $d$. If $h(t) = h_{0} + h_{1} \, t + \cdots +  h_{d} \, t^{d}$ is a nonzero polynomial  of degree at most $d$ with nonnegative integer coefficients and $h_{0} = 1$, then
\[
\frac{h(t)}{(1 - t)^{d + 1}} = \sum_{m \geq 0} g(m) t^{m},
\]
where
$g(m) = \sum_{i = 0}^{d} h_{i} \binom{m + d - i}{d}$ is a polynomial of degree $d$ with rational coefficients.
For every positive integer $n$, define $\U_{n}h(t)$ to be the polynomial of degree at most $d$ with integer coefficients satisfying
\[ \frac{\U_{n}h(t)}{(1- t)^{d + 1}} =
\sum_{m \geq 0} g(nm) \, t^{m},
\]
and write
$\U_{n}h(t) = h_{0}(n) + h_{1}(n) \, t + \cdots +  h_{d}(n) \, t^{d}$.
The (Hecke) operator $\U_{n}$ was studied by Gil and Robins in a more general setting \cite{GRHecke} and more recently by
Brenti and Welker \cite{BWVeronese}.
%Brenti and Welker showed that for $n$ sufficiently large, the coefficients  $\U_{n}h(t)$ are well-behaved in a sense to be defined.
The goal of this paper is to show that there exists a positive integer $n_{d}$, depending only on $d$, such that
$\U_{n}h(t)$ is well-behaved, in a sense to be defined, for $n \geq n_{d}$.

Our main motivating example comes from the theory of lattice point enumeration of polytopes.
More specifically, let $N$ be a lattice of rank $n$ and set $N_{\mathbb{R}} := N \otimes_{\mathbb{Z}} \mathbb{R}$.
A \emph{lattice polytope} $P \subset N_\R$ is the convex hull of finitely many points in $N$.
Fix a $d$-dimensional lattice polytope $P \subset N_\R$ and, for each positive integer $m$, let $f_{P}(m) := \# \left( mP \cap N \right)$ denote the
number of lattice points in the $m$'th dilate of $P$. A famous theorem of Ehrhart \cite{ehrhartpolynomial} asserts that
$f_{P}(m)$ is a polynomial in $m$ of degree $d$, called the
\emph{Ehrhart polynomial} of $P$, and $f_{P}(0) = 1$.
Equivalently, the generating series of $f_{P}(m)$ can be written in the form
\begin{equation*}
  \frac{ \delta_{P}(t) }{ (1 - t)^{d + 1} } =
\sum_{m \geq 0} f_{P}(m) \, t^{m}  \, ,
\end{equation*}
where $\delta_{P}(t) = \delta_{0} + \delta_{1}t + \cdots +  \delta_{d}t^{d}$ is a polynomial of degree at most $d$ with
integer coefficients,
called the \emph{$\delta$-polynomial} of $P$, and $\delta_{0} = 1$. We call $(\delta_{0}, \delta_{1}, \ldots, \delta_{d})$ the \emph{(Ehrhart) $\delta$-vector} of $P$; alternative names in the literature include \emph{Ehrhart $h$-vector} and \emph{$h^*$-vector} of $P$. Stanley proved that
the coefficients $\delta_{i}$ are nonnegative \cite{StaDecompositions}.
In this case,
$\U_{n} \delta_{P}(t) = \delta_{nP}(t)$ and we write
$\delta_{nP}(t) = \delta_{0}(n) + \delta_{1}(n) \, t + \cdots +  \delta_{d}(n) \, t^{d}$.

More generally, let $R = \oplus_{i \geq 0} R_{i}$ be a graded ring of dimension $d + 1$ and assume that $R_{0} = k$ is a field and that $R$ is finitely generated
over $R_{0}$. The $n$'th \emph{Veronese subring} of $R$ is the graded ring $R^{\langle n \rangle} = \oplus_{i \geq 0} R_{in}$. The behaviour of Veronese subrings for large $n$ has been studied by Backelin \cite{BacRates} and Eisenbud, Reeves and Totaro \cite{ERTInitial}.
The \emph{Hilbert function} of $R$ is defined by $H(R,m) = \dim_{k} R_{m}$, for each nonnegative integer $m$, and by a theorem of Hilbert \cite[Theorem 4.1.3]{BHCohen}, $H(R,m)$ is a polynomial in $m$ of degree $d$ for $m$ sufficiently large.
In fact, $H(R,m)$ is a polynomial for $m > a(R)$, where $a(R)$ is the \emph{$a$-invariant} of $R$ and is defined in terms of the local cohomology of $R$
\cite[Section 3.6]{BHCohen}.
 Observe that this implies that $H(R^{\langle n \rangle},m)$ is a polynomial in $m$ of degree $d$ for $n > a(R)$.
Assume that $R$ is Cohen--Macauley and that $R$ is a finite module over the $k$-subalgebra of $R$ generated by $R_{1}$.
If $H(R,m)$ is a polynomial in $m$ then it can be seen as in \cite[Corollary 4.1.10]{BHCohen} that
\begin{equation*}
  \frac{ h_{0} + h_{1}t + \cdots + h_{d}t^{d} }{ (1 - t)^{d + 1} } =
\sum_{m \geq 0} H(R,m) \, t^{m}  \, ,
\end{equation*}
for some nonnegative integers $h_{i}$, with $h_{0} = 1$. For every positive integer $n$, the numerator of the generating series of $H(R^{\langle n \rangle},m)$ has the form
$(1 - t)^{d + 1} \sum_{m \geq 0} H(R^{\langle n \rangle},m) \, t^{m} =
   \U_{n}(h_{0} + h_{1}t + \cdots + h_{d}t^{d})$.
   Returning to our previous example, if $N' = N \times \mathbb{Z}$ and $\sigma$ denotes the cone over $P \times \{ 1 \}$ in
$N'_{\mathbb{R}}$, then the semigroup algebra $R = k[ \sigma \cap N' ]$ is graded
by the projection $u: N' \rightarrow \mathbb{Z}$ and satisfies the above assumptions
\cite[Theorem 6.3.5]{BHCohen}. In this case, $H(R,m) = f_{P}(m)$ is the Ehrhart polynomial of $P$
and $\U_{n} (h_{0} + h_{1}t + \cdots + h_{d}t^{d}) = \delta_{nP}(t)$. %we recover our previous example.
%\comment{this was confusing to me... $P$ seems to appear out of nowhere here.}

A sequence of positive integers $(a_{0}, \ldots, a_{d})$ is
\emph{strictly log concave} if $a_{i}^{2} > a_{i - 1}a_{i + 1}$ for $1 \leq i \leq d - 1$ and is \emph{strictly unimodal}
if $a_{0} < a_{1} < \cdots < a_{j}$ and $a_{j + 1} > a_{j + 2} > \cdots > a_{d}$ for some $0 \leq j \leq d$.
One easily verifies that if
$(a_{0}, \ldots, a_{d})$ is strictly log concave then it is strictly unimodal.
An induction argument implies that if the polynomial $a_{0} + a_{1}t + \cdots + a_{d}t^{d}$ has negative real roots then
the sequence $(a_{0}, \ldots, a_{d})$ is
strictly log concave and hence strictly unimodal.
Brenti and Welker recently proved the following theorem  \cite[Theorem~1.4]{BWVeronese}.

\begin{theorem}[Brenti--Welker] \label{brentiwelthm}
For any positive integer $d$, there exist real numbers $\alpha_{1} < \alpha_{2} < \cdots < \alpha_{d - 1} < \alpha_{d} = 0$ such that,
if $h(t) = h_{0} + h_{1} \, t + \cdots +  h_{d} \, t^{d}$ is a polynomial of degree at most $d$ with nonnegative integer coefficients
and $h_{0} = 1$, then for $n$ sufficiently large, $\U_{n} h(t)$ has negative real roots
$\beta_{1}(n) < \beta_{2}(n) < \cdots < \beta_{d - 1}(n) < \beta_{d}(n) < 0$ and $\beta_{i}(n) \rightarrow \alpha_{i}$ as $n \rightarrow \infty$\footnote{In \cite{BWVeronese}, the notation $\beta_{i}(n)$ is used for the reciprocals of the roots.}.
\end{theorem}

Let $w = (w_{1},\ldots,w_{d})$ be a permutation of $d$
elements. A \emph{descent} of $w$ is an index $1 \le j
\le d- 1$ such that $w_{j + 1} < w_{j}$. If $A(d,i)$
denotes the number of permutations of $d$ elements with
$i -1$ descents, then the polynomial $A_{d}(t) =
\sum_{i = 1}^{d} A(d,i) \, t^{i}$ is called an
\emph{Eulerian polynomial} and the roots of
$\frac{A_{d}(t)}{t}$ are simple, real and strictly
negative \cite[p.~292, Exercise 3]{ComAdvanced}. We are
ready to state our main result and emphasise that the
real content of Theorem \ref{Lebron} is the statement that
the constants $m_d$ and $n_d$ below only depend on $d$ (and \emph{not} on $h(t)$).

\begin{theorem}\label{Lebron}
Fix a positive integer $d$ and let $\rho_{1} < \rho_{2} < \cdots < \rho_{d} = 0$ denote the roots of the Eulerian polynomial $A_{d}(t)$. There exist positive integers $m_{d}$ and $n_{d}$ such that, if $h(t)$ is a polynomial of degree  at most $d$
with nonnegative integer coefficients and $h_{0} = 1$, then for $n \geq n_{d}$,
$\U_{n} h(t)$ has negative real roots
$\beta_{1}(n) < \beta_{2}(n) < \cdots < \beta_{d - 1}(n) < \beta_{d}(n) < 0$ with $\beta_{i}(n) \rightarrow \rho_{i}$ as $n \rightarrow \infty$, and
the coefficients of $\U_{n}h(t)$
%the sequence
%$\{ h_{0}(n), h_{1}(n), \ldots, h_{d}(n) \}$
are positive, strictly log concave, % strictly unimodal,
and satisfy $h_{i}(n) < m_{d}h_{d}(n)$ for $0 \leq i \leq d$.
Furthermore, we may choose $n_{d}$ such that, if additionally
%\[
$h_{0} + \cdots + h_{i+1} \geq h_{d} + \cdots + h_{d - i}$ %\quad \textrm{ for } \quad 0 \le i \le \left\lfloor \tfrac d 2 \right\rfloor , \]
for $0 \le i \le \left\lfloor \tfrac d 2 \right\rfloor - 1$, then
%with at least one of these inequalities strict, then for $n \geq n_{d}$,
%$h_{i + 1}(n) > h_{d - i}(n)$ for $n \geq n_{d}$ and $0 \leq i \leq \left\lfloor \frac d 2 \right\rfloor - 1$.
% and satisfying
%\[
%h_{0} + \cdots + h_{i} \geq h_{d + 1} + \cdots + h_{d + 1 - i} \quad \textrm{ for } \quad 0 \le i \le \left\lfloor \tfrac d 2 %\right\rfloor , \]
%with at least one of the above inequalities strict, then for any $n \ge n_{d}$,
\begin{equation*}
h_{0} = h_{0}(n) < h_{d}(n) < h_{1}(n) < \cdots < h_{i}(n) < h_{d - i}(n) < h_{i + 1}(n) < \cdots <
h_{\lfloor \frac{d + 1}{2} \rfloor}(n) < m_{d} \, h_{d}(n) \, .
\end{equation*}
%In particular, the coefficents of $\U_{n}h(t)$  are strictly unimodal for $n \geq n_{d}$.
%if $P$ is a $d$-dimensional lattice polytope and $n \geq n_{d}$, then
%\begin{equation*}
%\delta_{0}(n) < \delta_{1}(n) < \delta_{2}(n) < \dots < \delta_{\lfloor \frac{ d+1 }{ 2 } \rfloor}(n) < m_{d} \, \delta_{d}(n) \, ,
%\end{equation*}
%\begin{equation*}
%\delta_{d}(n)  < \delta_{d-1}(n)  < \dots < \delta_{\lfloor \frac d 2 \rfloor + 1}(n) < m_{d} \, \delta_{d}(n) \, .
%\end{equation*}
\end{theorem}

If $h(t) = \delta_{P}(t)$ then assumptions of the above theorem hold by a result of Hibi \cite{HibSome},
%the inequalities (\ref{007}) were established by
%Hibi in ,
and we deduce the following corollary.

\begin{cor}\label{mainthm}
Fix a positive integer $d$ and let $\rho_{1} < \rho_{2} < \cdots < \rho_{d} = 0$ denote the roots of the Eulerian polynomial $A_{d}(t)$. There exists positive integers $m_{d}$ and $n_{d}$ such that,  if $P$ is a $d$-dimensional lattice polytope and $n \geq n_{d}$, then
$\delta_{nP}(t)$ has negative real roots
$\beta_{1}(n) < \beta_{2}(n) < \cdots < \beta_{d - 1}(n) < \beta_{d}(n) < 0$ with $\beta_{i}(n) \rightarrow \rho_{i}$ as $n \rightarrow \infty$, and
the coefficients of $\delta_{nP}(t)$
%the sequence
%$\{ h_{0}(n), h_{1}(n), \ldots, h_{d}(n) \}$
are positive, strictly log concave, % strictly unimodal 
and satisfy %$\delta_{i}(n) < m_{d}h_{d}(n)$ for $0 \leq i \leq d$.
%Furthermore, if
%\[
%$h_{0} + \cdots + h_{i} \geq h_{d} + \cdots + h_{d + 1 - i}$ %\quad \textrm{ for } \quad 0 \le i \le \left\lfloor \tfrac d 2 \right\rfloor , \]
%for $0 \le i \le \left\lfloor \tfrac d 2 \right\rfloor$,
%with at least one of the above inequalities strict, then we may choose $n_{d}$ such that, for $n \geq n_{d}$,
%$h_{i + 1}(n) > h_{d - i}(n)$ for $n \geq n_{d}$ and $0 \leq i \leq \left\lfloor \frac d 2 \right\rfloor - 1$.
% and satisfying
%\[
%h_{0} + \cdots + h_{i} \geq h_{d + 1} + \cdots + h_{d + 1 - i} \quad \textrm{ for } \quad 0 \le i \le \left\lfloor \tfrac d 2 %\right\rfloor , \]
%with at least one of the above inequalities strict, then for any $n \ge n_{d}$,
\begin{equation*}
1 = \delta_{0}(n) < \delta_{d}(n) < \delta_{1}(n) < \cdots < \delta_{i}(n) < \delta_{d - i}(n) < \delta_{i + 1}(n) < \cdots <
\delta_{\lfloor \frac{d + 1}{2} \rfloor}(n) < m_{d} \, \delta_{d}(n) \, .
\end{equation*}

%Fix a positive integer $d$. There exists positive integers $m_{d}$ and $n_{d}$ such that, if $P$ is a $d$-dimensional lattice polytope and $n \geq n_{d}$, then the sequence $\{ \delta_{1}(n), \ldots, \delta_{d}(n) \}$ is strictly log concave and
%\begin{equation*}
%1 = \delta_{0}(n) < \delta_{d}(n) < \delta_{1}(n) < \cdots < \delta_{i}(n) < \delta_{d - i}(n) < \delta_{i + 1}(n) < \cdots <
%\delta_{\lfloor \frac{d + 1}{2} \rfloor}(n) < m_{d} \, \delta_{d}(n) \, .
%\end{equation*}
%In particular, the Ehrhart $\delta$-vector of $nP$ is strictly unimodal for $n \geq n_{d}$.
\end{cor}

We also have the following application to Veronese subrings of graded rings.

\begin{cor}\label{veronesecor}
Fix a positive integer $d$ and let $\rho_{1} < \rho_{2} < \cdots < \rho_{d} = 0$ denote the roots of the Eulerian polynomial $A_{d}(t)$. There exists positive integers $m_{d}$ and $n_{d}$ such that,
if $R = \oplus_{i \geq 0} R_{i}$ is a finitely generated graded ring over a field $R_{0} = k$, which is
 Cohen--Macauley and module finite over the $k$-subalgebra of $R$ generated by $R_{1}$, and if  the Hilbert function
$H(R,m)$ is a polynomial in $m$ and we write
\begin{equation*}
  \frac{ \U_{n} h(t) }{ (1 - t)^{d + 1} } =
\sum_{m \geq 0} H(R^{\langle n \rangle},m) \, t^{m}  \, ,
\end{equation*}
for each positive integer $n$,
then for $n \geq n_{d}$,
$\U_{n} h(t)$ has negative real roots
$\beta_{1}(n) < \beta_{2}(n) < \cdots < \beta_{d - 1}(n) < \beta_{d}(n) < 0$ with $\beta_{i}(n) \rightarrow \rho_{i}$ as $n \rightarrow \infty$, and
the coefficients of $\U_{n}h(t)$
%the sequence
%$\{ h_{0}(n), h_{1}(n), \ldots, h_{d}(n) \}$
are positive, strictly log concave, % strictly unimodal 
and satisfy $h_{i}(n) < m_{d}h_{d}(n)$ for $0 \leq i \leq d$.
\end{cor}

It is an open problem to determine the optimal choices for the integers $m_d$ and $n_d$ in Theorem \ref{Lebron} and Corollaries \ref{mainthm} and \ref{veronesecor}.
In this direction, we show that for any positive integer $d$ and $n \geq d$, if $h(t)$ satisfies certain inequalities,  %the assumptions of Theorem \ref{ole},
then $h_{i + 1}(n) > h_{d - i}(n)$ for
$i = 0, \ldots, \left\lfloor \frac d 2 \right\rfloor - 1$ (Theorem \ref{Beijing}). In particular, this holds when $h(t) = \delta_{P}(t)$
(Example \ref{bigshot}).
%$\delta_{i + 1}(n) > \delta_{d - i}(n)$, hold in fact for $n \ge d$.

We now explain our original motivation for this paper.
A triangulation $\tau$ of the polytope $P$ with vertices in $N$ is \emph{unimodular} if for any simplex of $\tau$ with vertices $v_{0}, v_{1}, \ldots, v_{d}$, the vectors $v_{1} - v_{0}, \ldots, v_{d} - v_{0}$ form a basis of $N$.
%\item Whenever $G$ is a face of $\tau$, the cone over $G \times \{1\}$ in $(N \times \mathbb{Z})_{\mathbb{R}}$ is nonsingular.
%\end{enumerate}
%Equivalently,
%let $\sigma$ be the cone over $P \times \{1\}$ in $(N \times \mathbb{Z})_{\mathbb{R}}$.
%Let $\triangle$ be the fan over $\tau$, which refines $\sigma$. Then $\tau$ is unimodular if and only if $\triangle$ is nonsingular.
While every lattice polytope can be triangulated into lattice simplices, it is far from true that every lattice polytope admits a \emph{unimodular} triangulation (for an easy example, consider the convex hull of $(1,0,0)$, $(0,1,0)$, $(0,0,1)$, and $(1,1,1)$). The following theorem, however, says that we can obtain a unimodular triangulation if we allow our polytope to be dilated.

\begin{theorem}[Knudsen--Mumford--Waterman \cite{KKFMSDToroidal}]\label{oldies}
For every lattice polytope $P$, there exists an integer $n$ such that $nP$ admits a regular unimodular triangulation.
\end{theorem}

\noindent
For a general reference on triangulations, including regular ones, 
see \cite{leetriangulations}.
%Note that
If $P$ admits a unimodular triangulation, then
every multiple $nP$ admits such a triangulation
%,
%because any multiple of a unimodular simplex admits a
% unimodular triangulation
(this follows from the
general theory of Knudsen--Mumford triangulations; see
\cite[Remark 3.19]{brunsgub}). Thus Theorem
\ref{oldies} implies that $knP$ admits a
unimodular triangulation for $k \in \Z_{>0}$. There are
several conjectured stronger versions of
Theorem~\ref{oldies} (see, for example,
\cite{BGUnimodular,BGNormal}):

\begin{conj}\label{oldieconjecture}
\begin{enumerate}[{\rm (a)}]
\item For every lattice polytope $P$, there exists an integer $m$ such that $nP$ admits a regular unimodular triangulation for $n \ge m$.
\item For every $d \in \Z_{ >0 } $, there exists an integer $n_d$ such that, if $P$ is a $d$-dimensional lattice polytope, then $n_d P$ admits a regular unimodular triangulation.
\item For every $d \in \Z_{ >0 } $, there exists an integer $n_d$ such that, if $P$ is a $d$-dimensional lattice polytope, then $n P$ admits a regular unimodular triangulation for $n \ge n_d$.
\end{enumerate}
\end{conj}

\noindent
When $d = 1$ or $2$, every lattice polytope has a unimodular triangulation.
For $d = 3$, Conjecture \ref{oldieconjecture}(b) holds with $n_{3}= 4$ \cite{KSPrimitive}.

Conjecture \ref{oldieconjecture} was the first motivation for our paper, and the following result \cite[Theorem 1.3]{AthVectors} was the second.
\begin{theorem}[Athanasiadis--Hibi--Stanley]\label{athanasiadisthm}
If a $d$-dimensional lattice polytope $P$ admits a regular unimodular triangulation, then the $\delta$-vector of $P$ satisfies
\begin{enumerate}[{\rm (a)}]
\item $\delta_{i + 1} \ge \delta_{ d-i }$ for $0 \le i \le \lfloor \frac{ d }{ 2 } \rfloor - 1 \, ,$
\item $\delta_{ \lfloor \frac{ d+1 }{ 2 } \rfloor } \ge \delta_{ \lfloor \frac{ d+1 }{ 2 } \rfloor + 1 } \ge \dots \ge \delta_{ d-1 } \ge \delta_d \, ,$
\item $\delta_i \le \binom{ \delta_1 + i -1 }{ i }$ for $0 \le i \le d$.
\end{enumerate}
In particular, if the $\delta$-vector of $P$ is symmetric and $P$ admits a regular unimodular triangulation, then the $\delta$-vector is unimodal.
\end{theorem}

%\comment{A: Should we say something later about the last inequality? M: Yes, please. I remember you saying that this follows sort of trivially here. Could you add a meaningful remark? A: I am not sure where to include it but the point is: the reason it is true is that $\delta_{i}(n)$ is a polynomial in $n$ of degree $d$ with positive leading coefficient, while the fact that $\delta_{1}(n)$ is a polynomial in $n$ of degree $d$ with positive leading coefficient means that the right hand side is a polynomial in $n$ of degree $id$ with positive leading coefficient; so the right hand side has much faster growth. To prove that the inequality holds for $n \geq n_{d}$, you run the same (root bounding) argument as for the main theorem; using the same inequalities of Betke and McMullen. M: Naturally. It's not clear where to fit it in... let's leave it for now.}

\noindent
In fact, the first inequality in the above theorem holds under the weaker assumption that the boundary of $P$ admits a
regular unimodular triangulation \cite[Theorem 2.20]{YoInequalities}.
There are (many) lattice polytopes for which some of the inequalities of Theorem \ref{athanasiadisthm} fail and one may hope to
use Theorem \ref{athanasiadisthm} to construct a counter-example to Conjecture \ref{oldieconjecture}.
However, a consequence of Corollary \ref{mainthm} and its proof is that this approach can not possibly work. More precisely, one can show that there exists a positive integer
$n_{d}$ such that  if $n \geq n_{d}$, then the inequalities in Theorem \ref{athanasiadisthm} hold for $nP$.
% one can not hope to use Theorem \ref{athanasiadisthm} in order to provide a counterexample to Conjecture~\ref{oldieconjecture}.
%However, Theorem \ref{oldies} and Conjecture \ref{oldieconjecture} suggest that the inequalities of Theorem \ref{athanasiadisthm} hold for the $\delta$-vectors of sufficiently large \emph{dilates} of $P$. This motivates the main object of our paper: Our goal is to study the $\delta$-vector of $nP$ as $n$ increases.
%Given a $d$-dimensional lattice polytope $P$, we let
%\begin{equation*}
%\delta_{nP}(t) = \delta_{0}(n) + \delta_{1}(n) \, t + \cdots +  \delta_{d}(n) \, t^{d} .
%\end{equation*}
%Recall that a vector $(a_0, a_1, \dots, a_d)$ is \emph{strictly unimodal} if there exists an index $k$ such that $a_0 < a_1 < \dots < a_k > a_{ k+1 } > \dots > a_d$.
%Our main result is as follows.

%\begin{theorem}\label{mainthm}
%Fix a positive integer $d$. There exists positive integers $m_{d}$ and $n_{d}$ such that, if $P$ is a $d$-dimensional lattice polytope %and $n \geq n_{d}$, then
%\begin{equation*}
%1 = \delta_{0}(n) < \delta_{d}(n) < \delta_{1}(n) < \cdots < \delta_{i}(n) < \delta_{d - i}(n) < \delta_{i + 1}(n) < \cdots <
%\delta_{\lfloor \frac{d + 1}{2} \rfloor}(n) < m_{d} \, \delta_{d}(n) \, .
%\end{equation*}
%In particular, the Ehrhart $\delta$-vector of $nP$ is strictly unimodal for $n \geq n_{d}$.
%\end{theorem}

We end the introduction with a brief outline of the contents of the paper.
%We will prove Theorem \ref{mainthm} in a slightly more general setting (Theorem \ref{ole} below) in Section \ref{heckesection}. It follows from
In Section \ref{inequalities}, we develop some
inequalities between the coefficients of polynomials with certain properties, %which we develop in Section \ref{inequalities}.
and we remark that Theorem \ref{runner} might be interesting in its own right---it asserts that we can bound roughly half the coefficients of an Ehrhart polynomial in terms of the dimension of $P$ and the surface area of $P$.
In Section \ref{heckesection}, we express $h_{i}(n)$ as a sum of Eulerian polynomials for $1 \leq i \leq d$ and use this description to  establish our main results. In Section \ref{lowerboundsection}, we consider %lower 
bounds for $n_{d}$ and prove the aforementioned Theorem \ref{Beijing}.
We conclude in Section \ref{openquestionsection} with a conjecture on Ehrhart $\delta$-vectors.
%show that some of the inequalities in Corollary \ref{mainthm}, namely, $\delta_{i + 1}(n) > \delta_{d - i}(n)$, hold in fact for $n \ge d$.
%\comment{A: We should think about this a little. Can we make an accurate conjecture at least? Lagarias and Ziegler has some bounds on the volume of a lattice polytope with an interior lattice point. The bounds depend on $d$ and the number of interior lattice points $\delta_{d} > 0$. Lagarias told me that the bounds are doubly exponential (I haven't looked at the paper though). Should we expect similar behaviour? Remember that the volume of $P$ is the $\frac{ 1 }{d!}$ times sum of the coefficients of the $\delta$-vector. What about optimal bounds for small $d$? M: Good point. I'll try to think about this... It might be worthwhile, at any rate, to talk about the Lagarias--Ziegler paper in this connection.}

%%%%%%%%%%%%%%%%%%%%%%%%%%%%%%%%%%%%%%%%%%%%%%%%%%%%%%%%%%%%%%%

\section{Inequalities between Coefficients of Polynomials}\label{inequalities}

%The goal of this section is to recall some inequalities of Betke and McMullen between the coefficients of the Ehrhart polynomial of a lattice polytope.

Our setup in this section will be slightly more general than the one in the introduction. We fix the following notation throughout the paper.
Let  $h(t) = h_{0} + h_{1}t + \cdots +  h_{d + 1} t^{d + 1}$ be a nonzero polynomial of degree at most $d + 1$ with integer coefficients,
and write
\begin{equation}\label{star}
h_{0} + \sum_{m \geq 1} g(m) \, t^{m} = \frac{ h(t) }{ (1 - t)^{d + 1} } \, ,
\end{equation}
where
$g(m) = \sum_{i = 0}^{d + 1} h_{i} \binom{m + d - i}{d}$ is a polynomial with rational coefficients.
We write $g(t) = g_{d}t^{d} + g_{d -1}t^{d -1} + \cdots + g_{0}$ and will assume that $g_{d} = \frac{\sum_{i = 0}^{d + 1} h_{i}}{d!}$ is positive and hence bounded below by $\frac{1}{d!}$. One can verify that $h_{d + 1} = (-1)^{d}(g(0) - h_{0})$ and we will often assume that $h(t)$ is a polynomial of degree at most $d$, in which case $g(0) = h_{0}$.

%Consider a function $g: \mathbb{Z}_{> 0} \rightarrow \mathbb{Z}_{>0}$. We will assume that $g(m)$ is a polynomial in $m$ of degree $d$ with rational coefficients and write $g(m) = g_{d}t^{d} + g_{d -1}t^{d -1} + \cdots + g_{0}$. Observe that $g(m)$ has a positive leading coefficient.
%Equivalently, for any fixed integer $\alpha$, the generating series of $g(m)$ can be written in the form
%\begin{equation}\label{star}
%\alpha + \sum_{m \geq 1} g(m) \, t^{m} = \frac{ h(t) }{ (1 - t)^{d + 1} } \, ,
%\end{equation}
%where $h(t) = h_{0} + h_{1}t + \cdots +  h_{d + 1} t^{d + 1}$ is a polynomial of degree at most $d + 1$ with
%rational coefficients \cite[Chapter 4]{StaEnumerative}.
%One can verify that $h_{0} = \alpha$ and $h_{d + 1} = (-1)^{d}(g(0) - \alpha)$. We will often set $\alpha = g(0)$ so that $h(t)$ is a polynomial of degree at most $d$.
%Expanding (\ref{star}) gives
%\[
%\alpha + \sum_{m \geq 1} g(m) \, t^{m} = \frac{ h(t) }{ (1 - t)^{d + 1} } = \sum_{i = 0}^{d + 1} h_{i} t^{i} \sum_{m \geq 0}  \binom{d + m}{d}t^{m} = \sum_{i = 0}^{d + 1} h_{i}  \sum_{m \geq i}  \binom{d + m - i}{d}t^{m},
%\]
%and hence
%\begin{equation}\label{star2}
%g(m) = \sum_{i = 0}^{d + 1} h_{i} \binom{m + d - i}{d} \, .
%\end{equation}

\begin{exa}\label{strike}
A \emph{lattice complex} $K$ in a lattice $N$ is a simplicial complex in $N_{\mathbb{R}}$ whose vertices lie in $N$. A lattice complex is
\emph{pure} of dimension $r$ if all its maximal simplices have dimension $r$. Let $K$ be a pure lattice complex of dimension $r$ and, for each positive integer $m$, let
$f_{K}(m) := \# \left( mK \cap N \right) $
denote the
number of lattice points in the $m$'th dilate of $K$. Ehrhart's theorem implies that $f_{K}(m)$ is a polynomial in $m$ of degree $r$. If we write
%\begin{equation*}
$1 + \sum_{m \geq 1} f_{K}(m) \, t^{m} = \frac{ \delta_{K}(t) }{ (1 - t)^{r + 1} } \,$ ,
%\end{equation*}
%\comment{M: Shouldn't this $1$ on the left be replaced by the Euler characteristic of $K$? A: No, Betke and McMullen use $1$; they make a special mention of this in their paper.}
then Betke and McMullen \cite{BMLattice} showed that $\delta_{K}(t)$ has nonnegative coefficients if $K$ is homeomorphic to a ball or a sphere.
Moreover, $\delta_{K}(t)$ has degree at most $d$ when $K$ is homeomorphic to a ball and the coefficients of $\delta_{K}(t)$ are symmetric when $K$ is homeomorphic to a sphere.
%Also, the coefficients of $\delta_{K}(t)$ are nonnegative and symmetric when $K$ is homeomorphic to a $r$-sphere.
For example, a $d$-dimensional lattice polytope $P$ is homeomorphic to a $d$-ball and can be given the structure of a pure lattice complex of dimension $d$. Its boundary $\partial P$ is homeomorphic to a $(d-1)$-sphere and can be given the structure of a pure lattice complex of dimension $d -1$.
\end{exa}

The following inequalities and their proof are a slight generalisation of \cite[Theorem 6]{BMLattice}.
Recall that the \emph{Stirling number $S_{i}(d)$ of the first kind} is the coefficient of $t^{i}$ in $\prod_{j = 0}^{d - 1} (t - j)$; note that $(-1)^{ d-i } S_i(d) > 0$ for $i \ge 1$.

\begin{theorem}[Betke--McMullen]\label{inequality}
With the notation of \eqref{star}, if
$h_{i} \geq 0$  for  $0 \le i \le d + 1$, then for any $1 \le r \le d - 1$,
\[
g_{r} \leq (-1)^{d - r} S_{r}(d) \, g_{d} + \frac{(-1)^{d - r- 1} \, h_{0} \, S_{r + 1}(d)}{(d - 1)!} \, .
\]
\end{theorem}
\begin{proof}
%If $e(t) = \sum_{m \geq 0} e_{m} t^{m}$ and $f(t) = \sum_{m \geq 0} f_{m} t^{m}$ are power series in $t$ with integer coefficients, then we will write $e(t) \leq f(t)$ if $e_{m} \leq f_{m}$ for all $m \geq 0$.
%Using (\ref{star}) and the nonnegativity of the $h_{i}$,
%\[
%\alpha + \sum_{m \geq 1} g(m) \, t^{m} = \frac{ h(t) }{ (1 - t)^{d + 1} } = \sum_{i = 0}^{d + 1} h_{i} t^{i} \sum_{m \geq 0}  \binom{d + m}{d}t^{m} \leq h_{0}\sum_{m \geq 0}  \binom{d + m}{d}t^{m} + \sum_{i = 1}^{d + 1} h_{1}  \sum_{m \geq 1}  \binom{d + m - 1}{d}t^{m}.
%\]
%Comparing coefficients of $m^{r}$ in both sides of (\ref{star2}) yields
By definition,
$g_{r} = \sum_{i = 0}^{d + 1} h_{i} \binom{m + d - i}{d}_{r}$,
where $\binom{m + d - i}{d}_{r}$ denotes the coefficient of $m^{r}$ in $\binom{m + d - i}{d}$. Observe that
%\[
$\binom{m + d}{d}_{r} \geq \binom{m + d - 1}{d}_{r} \geq \binom{m + d - i}{d}_{r} \quad \textrm{ for } \quad 2 \leq i \leq d + 1$,
%\]
and hence, by the nonnegativity of the $h_{i}$,
%\[
$g_{r} \leq h_{0} \binom{m + d}{d}_{r} + \sum_{i = 1}^{d + 1} h_{i} \binom{m + d - 1}{d}_{r} \,$.
%\]
%By considering the leading term of both sides of (\ref{star2}), we see that
Using the fact that $d! \, g_{d} = \sum_{i = 0}^{d + 1} h_{i}$ and applying a binomial identity, we get
%\[
$g_{r} \leq h_{0} \binom{m + d - 1}{d - 1}_{r} + \, d! \, g_{d} \binom{m + d - 1}{d}_{r} \, $.
%\]
Observing that $\binom{m + d - 1}{d - 1}_{r}$ is the coefficient of $m^{r + 1}$ in $\frac{\prod_{j = 0}^{d - 1}(m +j)}{(d - 1)!}$, which
is the coefficient of $(-m)^{r + 1}$ in  $\frac{(-1)^{d}\prod_{j = 0}^{d - 1}(m - j)}{(d - 1)!}$, we conclude that
$\binom{m + d - 1}{d - 1}_{r} = \frac{(-1)^{d - r - 1}S_{r + 1}(d)}{(d - 1)!}$. Similarly, one can verify that
$\binom{m + d - 1}{d}_{r} = \frac{(-1)^{d - r}S_{r}(d)}{d!}$ and the result follows.
\end{proof}

\begin{exa}\label{night}
If $P$ is a $d$-dimensional lattice polytope, denote its Ehrhart polynomial by
%\[
$f_{P}(m) = c_{d}m^{d} + c_{d - 1}m^{d - 1} + \cdots + c_{0} \,$.
%\]
Basic facts of Ehrhart theory (see, e.g., \cite{BRComputing}) imply that $c_{d}$ is the normalised volume of $P$ and $c_{d - 1}$ is half the normalised surface area of $P$. In this case, $h(t) = \delta_{P}(t)$ is the Ehrhart $\delta$-polynomial of $P$ and $\delta_{0} = 1$.
Since the coefficients of $\delta_{P}(t)$ are nonnegative \cite{StaDecompositions}, Theorem \ref{inequality} implies that
the coefficients $c_{i}$ can be bounded in terms of $d$ and the volume of $P$ (a fact that follows also, e.g., from \cite{lagariasziegler}).
\end{exa}

We can strengthen these inequalities if we put further restrictions on the coefficients $h_{i}$. We will need the following lemmas, the first of which is motivated by similar results in~\cite{YoInequalities}. 
%\footnote{In \cite{YoInequalities}, the author considered a decomposition} 

\begin{lemma}\label{banana}
A polynomial $h(t) = h_{0} + h_{1}t + \cdots + h_{d + 1}t^{d + 1}$ with integer coefficients
has a unique decomposition
%\begin{equation*}
$h(t) = a(t) + b(t)$,
%\end{equation*}
where $a(t)$ and $b(t)$  are polynomials with integer coefficients satisfying
$a(t) = t^{d} \, a(\frac 1 t)$ and $b(t) = t^{d + 1} \, b(\frac 1 t)$.
%Moreover, if $a_{i}$ denotes the coefficient of $t^{i}$ in $a(t )$, then
%\[
%a_{i} = h_{0} + h_{1} + \cdots + h_{i} - h_{d + 1} - h_{d} - \cdots - h_{d + 1 - i}.
%\]
\end{lemma}
\begin{proof}
Let $a_{i}$ and $b_{i}$ denote the coefficients of $t^{i}$ in $a(t)$ and $b(t)$ respectively, and set
\begin{equation}\label{coeff1}
a_{i} = h_{0} + \cdots + h_{i} - h_{d + 1} - \cdots - h_{d + 1 - i} \, ,
\end{equation}
\begin{equation*}%\label{coeff2}
b_{i} = -h_{0} - \cdots - h_{i - 1} + h_{d + 1} + \cdots + h_{d + 1 - i} \, .
\end{equation*}
We see that $h(t) = a(t) + b(t)$
and %compute %using (\ref{pinochet}),
%a_{i} + b_{i - 1} &= h_{0} + \cdots + h_{i} - h_{d} - \cdots - h_{d - i + 1}
%-h_{0} - \cdots - h_{i - 1} + h_{d} + \cdots + h_{d - i+ 1} \\
%&= h_{i} \, , \\
\[a_{i} - a_{d - i} = h_{0} + \cdots + h_{i} - h_{d + 1} - \cdots - h_{d - i + 1}
- h_{0} -  \cdots - h_{d - i} + h_{d + 1} + \cdots + h_{i + 1} = 0,\]
\[b_{i} - b_{d + 1 - i} = -h_{0} - \cdots - h_{i - 1} + h_{d + 1} + \cdots + h_{d + 1 - i}
+ h_{0} + \cdots + h_{d - i} - h_{d + 1} - \cdots - h_{i} = 0,
\]
for $0 \le i \le d + 1$. Hence we obtain our desired decomposition and one easily verifies the uniqueness assertion.
\end{proof}

\begin{rem}\label{geo}
Alternatively, to prove the above lemma, one can check that
%\[
$a(t) = \frac{h(t) - t^{d + 1}h(t^{-1})}{1 -t}$ %\qquad \text{ and } \qquad b(t) = \frac{-t \, h(t) + t^{d + 1}h(t^{-1})}{1 -t} \, .
%\]
and $b(t) = \frac{-t \, h(t) + t^{d + 1}h(t^{-1})}{1 -t} \, $.
\end{rem}

\begin{rem}\label{Memphis}
It follows from \eqref{coeff1} that $a(t)$ is nonzero with nonnegative integer coefficients
%the coefficients of $a(t)$ are positive
if and only if
%\[
$h_{0} + \cdots + h_{i} \leq h_{d + 1} + \cdots + h_{d + 1 - i}$ %\quad \textrm{ for } \quad 0 \le i \le \left\lfloor \tfrac d 2 \right\rfloor , \]
for $0 \le i \le \left\lfloor \tfrac d 2 \right\rfloor$,
with at least one of these inequalities strict. The coefficients of $a(t)$ are positive if and only if each of the above inequalities are strict.
%Also, $a(t)$ has degree $d$ and positive leading coefficient if and only if $h_{0} > h_{d + 1}$.
Since
$a_{i + 1} - a_{i} = h_{i + 1} - h_{d - i}$, we see that
the coefficients of $a(t)$ are unimodal (resp.\ strictly unimodal) if and only if $h_{i + 1} \ge h_{d - i}$ (resp.\ $h_{i + 1} > h_{d - i}$)
for $0 \leq i \leq \left\lfloor \tfrac d 2 \right\rfloor - 1$.
\end{rem}

\begin{exa}\label{JayHawk}
If $P$ is a $d$-dimensional lattice polytope and we write $\delta_{P}(t) = a(t) + b(t)$ as in Lemma \ref{banana}, then
%, after comparing (\ref{coeff1}) above with $(17)$ in \cite{YoInequalities},
\cite[Theorem 2.14]{YoInequalities} implies that
%\begin{equation*} %\label{lanquin}
$1 = a_{0} \leq a_{1} \leq a_{i}$ % \quad \textrm{ for } \quad 2 \le i \le d - 1 .
%\end{equation*}
for  $2 \le i \le d - 1$.
In particular, $a(t)$ has degree $d$ and positive integer coefficients.
We have $\delta_{0} = 1 > \delta_{d + 1} = 0$ and,
by the above remark,
%\[
$\delta_{0} + \cdots + \delta_{i} \leq \delta_{d} + \cdots + \delta_{d + 1 - i}$, %\quad \textrm{ for } \quad 1 \le i \le \left\lfloor \tfrac d 2 \right\rfloor , \]
for $1 \le i \le \left\lfloor \tfrac d 2 \right\rfloor$.
The latter inequalities were proved by Hibi \cite{HibSome}, and the fact that the coefficients of $a(t)$ are positive implies that all of the inequalities are in fact strict.
\end{exa}

\begin{lemma}\label{dinner}
With the notations of \eqref{star} and Lemma \ref{banana}, if $g'(m) = g(m) - (-1)^{d}g(-m)$, then
\[
g(0) + \sum_{m \ge 1} g'(m) = \frac{ a(t)}{(1 - t)^{d} } \, .
\]
\end{lemma}
\begin{proof}
It is a standard fact (see, e.g., \cite[Exercise 4.6]{BRComputing}) that if
$G(t) = \sum_{m \geq 0} g(m) \, t^{m}$, then
$\sum_{m \geq 1} g(-m) \, t^{m} = -G(t^{-1})$.  Since $G(t) = \frac{h(t)}{(1 - t)^{d + 1}} + g(0) - h_{0}$, we see that
$\sum_{m \geq 1} g(-m) \, t^{m} = h_{0} - g(0) - \frac{(-t)^{d + 1}h(t^{-1})}{(1 - t)^{d + 1}}$.
We compute, using Remark \ref{geo},
%\[
$g(0) + \sum_{m \ge 1} g'(m) = \frac{h(t)}{(1 - t)^{d + 1}} + (-1)^{d} \, \frac{(-t)^{d + 1}h(t^{-1})}{(1 - t)^{d + 1}} = \frac{ a(t)}{(1 - t)^{d} } \, $. \qedhere
%\]
\end{proof}

\begin{lemma}\label{deadset}
With the notations of \eqref{star} and Lemma \ref{banana}, %if $h_{i} \geq 0$  for  $0 \le i \le d + 1$ then
$g_{d} \geq \frac{1}{d!}$. Furthermore, if
%\[
$h_{0} + \cdots + h_{i} \geq h_{d + 1} + \cdots + h_{d + 1 - i}$ %\quad \textrm{ for } \quad 0 \le i \le \left\lfloor \tfrac d 2 \right\rfloor , \]
for $0 \le i \le \left\lfloor \tfrac d 2 \right\rfloor$,
with at least one of these inequalities strict, then
%\[
%h_{0} + \cdots + h_{i} > h_{d + 1} + \cdots + h_{d + 1 - i} \textrm{ for } 0 \le i \le \lfloor d/2 \rfloor, \] then
$g_{d - 1} \ge \frac{1}{2(d - 1)!}$.
\end{lemma}
\begin{proof}
By assumption, $g_{d} = \frac{1}{d!}\sum_{i = 0}^{d + 1}h_{i}$ is positive and hence bounded below by $\frac{1}{d!}$.
%Comparing leading terms of both
%sides of (\ref{star2}) yields $g_{d} = \frac{1}{d!}\sum_{i = 0}^{d + 1}h_{i} \geq \frac{1}{d!}$, since $g_{d}$ is positive by assumption.
One verifies that the coefficient of $m^{d - 1}$ in $\binom{m + d - i}{d}$ is $\frac{d + 1 - 2i}{2(d - 1)!}$ and  hence
%, after comparing coefficients of $m^{d - 1}$ in both sides of (\ref{star2}), we see that
$g_{d - 1} = \sum_{i = 0}^{d + 1} h_{i}\frac{d + 1 - 2i}{2(d - 1)!}$.
By assumption, we have %\[
$h_{0} + \cdots + h_{i} \ge h_{d + 1} + \cdots + h_{d + 1 - i}$ % \quad \textrm{ for } \quad 0 \le i \le d, \]
for $0 \le i \le d$,
with at least one of these inequalities strict.
Summing these inequalities gives
$\sum_{i = 0}^{d + 1} (d + 1- i)h_{i} \geq \left( \sum_{i = 0}^{d + 1} i \, h_{i} \right) + 1 \, $
and we conclude that $g_{d - 1} = \sum_{i = 0}^{d + 1} h_{i}\frac{d + 1 - 2i}{2(d - 1)!} \geq \frac{1}{2(d - 1)!} \, .$
\end{proof}

%\begin{rem}\label{wiff}
%If we assume that
%\[
%h_{0} + \cdots + h_{i} > h_{d + 1} + \cdots + h_{d + 1 - i} \textrm{ for } 0 \le i \le \lfloor d/2 \rfloor, \]
%then the proof above shows that $g_{d - 1} \geq \frac{ d + 1}{2(d - 1)!}$.
%If $P$ is a $d$-dimensional lattice polytope and $h(t) = \delta_{P}(t)$, then the assumptions of the above lemma hold by Remark \ref{Memphis} and  Example \ref{JayHawk}, and, using Example \ref{night}, we recover the well-known fact that the normalised surface area of $P$ is at least $\frac{ d + 1}{(d - 1)!}$.
%\end{rem}

\begin{rem}
%If in the hypothesis of Lemma \ref{deadset} we only assume that
It follows from the proof of Lemma \ref{deadset} and \eqref{coeff1} that if $a(t)$ has degree $d$ and positive integer coefficients, then
%\[
%h_{0} + \cdots + h_{i} \ge h_{d + 1} + \cdots + h_{d + 1 - i} \textrm{ for } 0 \le i \le \lfloor d/2 \rfloor, \]
%and at least one of the above inequalities is strict, then $g_{d - 1} \geq \frac{1}{(d - 1)!}$.
$g_{d - 1} \ge \frac{d + 1}{2(d - 1)!}$. If $P$ is a $d$-dimensional lattice polytope and $h(t) = \delta_{P}(t)$, then, by Examples \ref{night} and \ref{JayHawk}, we recover the well-known fact that the normalised surface area of $P$ is at least $\frac{ d + 1}{(d - 1)!}$.
\end{rem}

In the case when $P$ is a $d$-dimensional lattice polytope and $h(t) = \delta_{P}(t)$, the existence of the following inequalities was
suggested by Betke and McMullen in \cite{BMLattice}.

\begin{theorem}\label{runner}
With the notation of \eqref{star}, if
$h_{0} + \cdots + h_{i} \geq h_{d + 1} + \cdots + h_{d + 1 - i}$ %\quad \textrm{ for } \quad 0 \le i \le \left\lfloor \tfrac d 2 \right\rfloor , \]
for $0 \le i \le \left\lfloor \tfrac d 2 \right\rfloor$,
with at least one of these inequalities strict, then
\[
g_{d - 1 - 2r}  \leq S_{d - 1 - 2r}(d - 1) \, g_{d - 1} - \frac{(h_{0} - h_{d + 1})S_{d - 2r}(d - 1)}{2(d - 2)!} \quad \textrm{ for } \quad 1 \le r \le
\left\lfloor \tfrac{d - 1}{2} \right\rfloor .
\]
\end{theorem}
\begin{proof}
By   Remark \ref{Memphis} and Lemma \ref{deadset},
the polynomial $g'(m) = g(m) - (-1)^{d}g(-m)$ has degree $d - 1$ and positive leading coefficient and, by Lemma \ref{dinner},
%\[
$g(0) + \sum_{m \ge 1} g'(m) \, t^m = \frac{ a(t)}{(1 - t)^{d} } \, $.
%\]
Since the coefficients of $a(t)$ are nonnegative by Remark \ref{Memphis}, applying Theorem \ref{inequality} to $a(t)$ yields inequalities on $g'(t)$ and hence on $g(t)$, namely,
%\[
$g_{d - 1 - 2r}  \leq S_{d - 1 - 2r}(d - 1) \, g_{d - 1} - \frac{a_{0} \, S_{d - 2r}(d - 1)}{2(d - 2)!}$ %\quad \textrm{ for } \quad 1 \le r \le
%\left\lfloor \tfrac{d - 1}{2} \right\rfloor ,
%\]
for $1 \le r \le
\left\lfloor \tfrac{d - 1}{2} \right\rfloor$ ,
where  $a_{0} = h_{0} - h_{d + 1}$ by~(\ref{coeff1}).
\end{proof}

\begin{exa}
If $P$ is a $d$-dimensional lattice polytope and $h(t) = \delta_{P}(t)$, then the assumptions of the above theorem hold by Remark \ref{Memphis} and  Example \ref{JayHawk}, and
hence we can bound the coefficients $c_{d - 1 - 2i}$ in terms of $d$ and the normalised surface area $2 \, c_{d -1}$ of $P$ (recalling Example \ref{night}). Betke and McMullen remark in \cite{BMLattice} that there are examples showing that a similar bound for $c_{d - 2}$ in terms of $d$ and $c_{d -1}$ does not exist.
\end{exa}

%%%%%%%%%%%%%%%%%%%%%%%%%%%%%%%%%%%%%%%%%%%%%%%%%%%%%%%%%%

\section{The Action of $\U_n$ on Integer Polynomials}\label{heckesection}

We will continue with the notation of the previous section and assume from now on that $h_{0} = 1$ and that $h(t)$ has degree at most $d$, so that
\begin{equation}\label{star3}
\sum_{m \geq 0} g(m) \, t^{m} = \frac{ h(t) }{ (1 - t)^{d + 1} } \, .
\end{equation}
%Recall that $g(m) = g_{d}m^{d} + g_{d -1}m^{d -1} + \cdots + g_{0}$ is a
%polynomial of degree $d$ with $g_{d} > 0$. We will assume from now on that $h_{0} = g(0)$, so that
%\begin{equation}\label{star3}
%\sum_{m \geq 0} g(m) \, t^{m} = \frac{ h(t) }{ (1 - t)^{d + 1} } \, ,
%\end{equation}
%where $h(t) = h_{0} + h_{1}t + \cdots +  h_{d} t^{d}$ is a polynomial of degree at most $d$ with
%integer coefficients. %Observe that $h_{0} = g(0)$.
Fix a positive integer $n$, and recall that $\U_{n}h(t)$ is the polynomial of degree at most $d$ with integer coefficients satisfying
%\[
$\sum_{m \geq 0} g(nm) \, t^{m} = \frac{\U_{n}h(t)}{(1- t)^{d + 1}} \, $.
%\]
%The (Hecke) operator $\U_{n}$ was studied by Brenti and Welker in \cite{BWVeronese} and in a more general setting by Gil and Robins in \cite{GRHecke}. 
We will write
%\[
$\U_{n}h(t) = h_{0}(n) + h_{1}(n) \, t + \cdots +  h_{d}(n) \, t^{d}$.
%\]
The goal of this section is to describe the behaviour of $\U_{n}h(t)$ for sufficiently large $n$.

\begin{exa}
If $P$ is a $d$-dimensional lattice polytope and we set $g(m) = f_{P}(m)$, then, with the notation of the introduction, $h(t) = \delta_{P}(t)$,
$\U_{n}h(t) = \delta_{nP}(t)$ and $h_{i}(n) = \delta_{i}(n)$.
\end{exa}

%shows that the coefficients
%of $\U_{n}h(t)$ are nonnegative.
The following well-known lemma should be compared with \cite[Theorem 1.1]{BWVeronese}.

\begin{lemma}\label{lions}
If $\E_{n}$ is the linear operator that takes a polynomial as input, discards its terms with powers
that are not divisible by $n$, and divides each remaining power by $n$, then
\[
\U_{n}h(t) =  \E_{n} \left( h(t) \, (1 + t + \cdots + t^{n -1})^{d + 1} \right).
\]
\end{lemma}
\begin{proof}
We extend $E_{n}$ to an operator on power series: given a degree-$d$ polynomial $h$, construct the polynomial $g$ such that $\sum_{ m \ge 0 } g(m) \, t^m = \frac{ h(t) }{ (1-t)^{ d+1 } }$. Applying $E_n$ to this rational generating function gives
\[ \sum_{m \geq 0} g(nm) \, t^{m} =
\E_{n} \left(\frac{h(t)}{(1 - t)^{d + 1} }\right) =
\E_{n} \left(\frac{h(t) \, (1 + t + \cdots + t^{n -1})^{d + 1}}{(1 - t^{n})^{d + 1} }\right) =
\frac{\E_{n} \left( h(t) \, (1 + t + \cdots + t^{n -1})^{d + 1} \right)}{(1 - t)^{d + 1}} \, . \qedhere
\]
\end{proof}

It follows from the definition that $h_{1}(n) = g(n) - (d + 1)$ is a polynomial in $n$ of degree $d$ with positive leading coefficient.
Our next goal will be to show that $h_{i}(n)$ is a polynomial in $n$ of degree $d$ with positive leading coefficient
for $1 \le i \le d$.
% in a lattice $N$.
%We will use the notation from the introduction.
%In particular, we write
%\begin{equation*}
%\delta_{P}(t) = \delta_{d}t^{d} + \delta_{d-1}t^{d-1} + \ldots + \delta_{0}.
%\end{equation*}
%We saw in Section \ref{bound} that $h_{0} = 1$, $h_{1} = |P \cap N| - (d + 1)$ and
%$h_{d}$ is the number of interior lattice points of $P$.
%We will investigate the behaviour of the $\delta$-vector of $nP$ as $n$ grows.
%Let $P$ be a $d$-dimensional lattice polytope.
%Recall that we write the $\delta$-polynomial of $nP$ as
%\begin{equation*}
%\delta_{nP}(t) = \delta_{d}(n) \, t^{d} + \delta_{d-1}(n) \, t^{d-1} + \cdots + \delta_{0}(n) \, ,
%\end{equation*}
%for each positive integer $n$. %From the above discussion and Ehrhart Reciprocity,
%Basic facts of Ehrhart theory (see, e.g., \cite{BRComputing}) imply that $\delta_{0}(n) \equiv 1$,  $\delta_{1}(n) = \# \left( nP \cap N \right) - (d+1) = f_{P}(n) - (d + 1)$, and $\delta_{d}(n) = \# \left( n P^\circ \cap N \right) = (-1)^{d}f_{P}(-n)$.
%(Here $P^\circ$ denotes the relative interior of $P$.)
%In particular, $\delta_0(n)$, $\delta_1(n)$, and $\delta_d(n)$ are polynomials in $n$. We will show that all $\delta_j(n)$ are polynomials, and for $1 \le j \le d$, they have degree $d$ and a positive leading coefficient.
%Let $w = (w_{1},\ldots,w_{d})$ be a permutation of $d$ elements.
%A \emph{descent} of $w$ is an index $1 \le j \le d- 1$ such that $w_{j + 1} < w_{j}$.
%Let $A(d,i)$ be the number of permutations of $d$ elements with $i -1$ descents. Note that $A(d,i)$ is zero unless $1 \le i \le d$.
Now we recall the Eulerian numbers $A(d,i)$ from the introduction; they
are positive and symmetric in the sense that $A(d, i) = A(d, d + 1 - i) \geq 1$ for $1 \leq i \leq d$ \cite[p. 242]{ComAdvanced}.
The nonzero roots of the Eulerian polynomial $A_{d}(t) = \sum_{i = 1}^{d} A(d,i)t^{i}$
are real and negative \cite[p.~292, Exercise 3]{ComAdvanced}, and consequently we have
%A short induction implies the well-known fact that if a polynomial has positive, rational coefficients and real, negative roots, then the coefficients of the polynomial are strictly log concave and hence strictly unimodal. That is,
\[
A(d,i)^{2} > A(d,i - 1)A(d,i + 1) \textrm{ for } 2 \leq i \leq d - 1
\]
\begin{equation*}
1 = A(d,1) < A(d,2) < \dots < A(d, \lfloor \tfrac{ d+1 }{ 2 } \rfloor) \, ,
\end{equation*}
\begin{equation*}
1 = A(d,d) < A(d,d - 1) < \dots < A(d, \lfloor \tfrac d 2 \rfloor + 1) \, .
\end{equation*}

If we set $g(m) = m^{d}$, then $h(t) = (1 - t)^{d + 1}\sum_{m \geq 0} m^{d} \, t^{m} = A_{d}(t)$ \cite[p.244]{ComAdvanced} and $\U_{n}A_{d}(t) = (1 - t)^{d + 1}\sum_{m \geq 0} (nm)^{d} \, t^{m}  = n^{d}A_{d}(t)$.
With the convention that $A_{0}(t) = 1$, we deduce the following lemma.

\begin{lemma}\label{Celtics}
If $g(m) = \sum_{j = 0}^{d} g_{j}m^{j}$  then
%$h(t) = \sum_{i = 0}^{d} g_{i}A_{i}(t)(1 - t)^{d - i}$
$\U_{n} h(t) = \sum_{j = 0}^{d} g_{j}A_{j}(t)(1 - t)^{d - j}n^{j}$, for every positive integer $n$. In particular, for $1 \leq i \leq d$,
$h_{i}(n)$ is a polynomial in $n$ of degree $d$ of the form
\begin{equation*}
h_{i}(n) = A(d,i) \, g_{d} \, n^{d} + (A(d-1,i) - A(d - 1, i - 1)) \, g_{d-1} \, n^{d-1} + O(n^{d -2}) \, .
\end{equation*}
\end{lemma}
\begin{proof}
We compute
\[
\U_{n} h(t) = (1 - t)^{d + 1}\sum_{m \geq 0} g(nm)t^{m} = (1 - t)^{d + 1}\sum_{j = 0}^{d} g_{j}n^{j} \sum_{m \geq 0} m^{j}t^{m}
= \sum_{j = 0}^{d} (1 - t)^{d - j}g_{j}n^{j} A_{j}(t),
\]
and the second statement follows.
\end{proof}

By Lemma \ref{Celtics} and the strict log concavity and strict unimodality of the Eulerian numbers, the integers $h_{i}(n)$ are strictly log concave and strictly unimodal for $n$ sufficiently large. Moreover, by the symmetry of the Eulerian numbers,
%\[
$h_{i + 1}(n) - h_{d - i}(n) = 2 (A(d-1,i) - A(d - 1, i - 1)) \, g_{d-1} \, n^{d-1} +  O(n^{d -2}) \, $.
%\]
Hence, if
$g_{d - 1} > 0$ then the strict unimodality of the Eulerian numbers implies that $h_{i + 1}(n) > h_{d - i}(n)$ for $n$ sufficiently
large and $0 \leq i \leq \left\lfloor \frac d 2 \right\rfloor - 1$.  In a similar direction, Brenti--Welker's Theorem \ref{brentiwelthm} says that for $n$ sufficiently large, $\U_{n} h(t)$ has negative real roots.
We will now consider the  existence of % lower 
bounds for such $n$.
%We will use the following inequality \cite[Theorem 6]{BMLattice}. It bounds the coefficients of the Ehrhart polynomial of $P$ in terms of the volume of $P$
%and its dimension $d$.
%\begin{theorem}[Betke--McMullen]\label{swannies}
%Let $P$ be a $d$-dimensional lattice polytope. Then for $1 \le j \le d-1$,
%\begin{equation*}
%c_{j} \leq (-1)^{d - j}s(d,j) \, c_{d} + (-1)^{d -j -1} \, \frac{ s(d,j+ 1) }{ (d -1)! } \, ,
%\end{equation*}
%where the $s(d,j)$ are the Stirling numbers of the first kind.
%\end{theorem}
We will use the following result of Cauchy (see, for example, \cite[Chapter VII]{MarGeometry}).

\begin{lemma}\label{mamacita}
Let $p(n) = p_{d} \, n^{d} + p_{d-1} \, n^{d -1 } + \dots + p_{0}$ be a polynomial of degree $d$ with real coefficients. The complex roots of $p(n)$ lie in the open disc
\begin{equation*}
\left\{ z \in \mathbb{C} : \, |z| < 1 + \max_{0 \leq j \leq d} \left| \frac{ p_{j} }{ p_{d} } \right| \right\} .
\end{equation*}
\end{lemma}

We are now ready to prove our main result.
Our method of proof should be compared with the proof of \cite[Theorem 1.2(a)]{BDDPSCoefficients}, which gives a bound
on the norm of the roots of the Ehrhart polynomial of a lattice polytope, and the proof of \cite[Lemma 4.7]{BWfVectors}.

\begin{proof}[Proof of Theorem \ref{Lebron}]
By Lemma \ref{Celtics}, $\U_{n} h(t) = n^{d}g_{d} \left( A_{d}(t) + \sum_{j = 0}^{d - 1} \frac{g_{j}}{g_{d}n^{d -j}}A_{j}(t)(1 - t)^{d - j} \right)$.
Fix $0 < \epsilon \ll 1$
%\comment{I think you mean $0 < \epsilon \ll 1$}
and let $M_{i} = \max_{t \in [\rho_{i}, \rho_{i} + \epsilon]} |A_{d}(t)|$
and  $M_{i}' = \max_{t \in [\rho_{i}  - \epsilon, \rho_{i}]} |A_{d}(t)|$ for $1 \le i \le d$.
 Since $g_{d} \geq \frac{ 1 }{ d! }$  by Lemma \ref{deadset}, and $|S_j(d)| \leq d!$, Theorem \ref{inequality} implies that
\begin{equation}\label{Hamilton}
\left| \frac{ g_{j} }{ g_{d} } \right| \le \left| (-1)^{d - j} S_j(d) + (-1)^{d -j -1} \frac{ S_{ j+1 } (d) } { g_{d}(d -1)! }
\right| \le d! + d! \, d \, ,
\end{equation}
for $1 \le j \le d - 1$. Hence, there exists a positive integer $N = N(d, \epsilon)$ such that  if $n \geq N(d, \epsilon)$, then
 \[ \max_{t \in [\rho_{i}, \rho_{i} + \epsilon]} \left| \sum_{j = 0}^{d - 1} \frac{g_{j}}{g_{d}n^{d -j}}A_{j}(t)(1 - t)^{d - j} \right| < M_{i}
\qquad \textrm{ and } \qquad \max_{t \in [\rho_{i} - \epsilon, \rho_{i}]} \left| \sum_{j = 0}^{d - 1} \frac{g_{j}}{g_{d}n^{d -j}}A_{j}(t)(1 - t)^{d - j} \right| < M_{i}' \, , \]
for $1 \le i \le d$. Since $\rho_i$ is a simple root of $A_d(t)$, it follows that for $1 \le i \le d$, there exists  $t_{i} \in [\rho_{i}, \rho_{i} + \epsilon]$ and
$t_{i}' \in [\rho_{i} - \epsilon, \rho_{i} ]$ such that  $\U_{n} h(t_{i}) \neq 0$ and $\U_{n}h(t_{i}') \neq 0$ have different signs.
Observe that since %$A_{j}(0) = 0$ for $1 \le j \le d$ and $A_{0}(0) = 1$,
$\U_{n} h(0) = h_{0} > 0$, we may and will set $t_{d} = \rho_{d} = 0$. We conclude that, for $n \geq N(d, \epsilon)$,
 we may choose $\beta_{i}(n) \in (t_{i}', t_{i})$ and
the first assertion follows. Note that if $\U_{n} h(t)$ has negative real roots, then it follows that
the coefficients of $\U_{n}h(t)$
are positive, strictly log concave and strictly unimodal.

By Lemma \ref{Celtics}, if we set $m_{d} = A(d, \lfloor \frac{ d+1 }{ 2 } \rfloor) + 1$, then
%the polynomial $m_{d} \, h_{d}(n) - h_{i}(n)$
\[
m_{d} \, h_{d}(n) - h_{i}(n) = (m_{d} - A(d,i))\, g_{d} d %\, n_{d} 
+ \sum_{j = 0}^{d -1 } \lambda_{j}(d)  \, g_{j} \, n^{j},
\]
for $0 \leq i \leq d$, where $\lambda_{j}(d)$ is a function of $d$, for $0 \le j \le d - 1$. By the strict unimodality of the Eulerian numbers,  $m_{d} \, h_{d}(n) - h_{i}(n)$ is a polynomial of degree $d$ with positive leading term.
It follows from Lemma \ref{mamacita} and (\ref{Hamilton}) that we can bound the absolute value of the roots of $m_{d} \, h_{d}(n) - h_{i}(n)$ in terms of $d$,
and we deduce the second assertion.
Similarly, we can bound the absolute values of the  roots of %$h_{i + 1}(n) - h_{i}(n)$ in terms of $d$ for $1 \le i \le \lfloor \frac{ d+1 }{ 2 } \rfloor$
%and bound the roots of $h_{i}(n) - h_{i + 1}(n)$ and
$h_{d - i}(n) - h_{i}(n)$ in terms of $d$ for $0 \le i \le \lfloor \frac{d - 1}{ 2} \rfloor$.
%For $i = 1, \dots, d$, the polynomial $m_{d} \, h_{d}(n) - h_{i}(n)$ has leading coefficient
%$(m_{d} - A(d,i))\, g_{d} \, n_{d}$. If we set $m_{d} = A(d, \lfloor \frac{ d+1 }{ 2 } \rfloor) + 1$, then it follows from the strict unimodality of the Eulerian numbers that $m_{d} \, \delta_{d}(n) - \delta_{i}(n)$ is a polynomial of degree $d$ with positive leading term. The
%above argument shows that we can bound the roots of $m_{d} \, \delta_{d}(n) - \delta_{i}(n)$ in terms of $d$.

Now assume that
%\[
%$h_{0} + \cdots + h_{i} \geq h_{d + 1} + \cdots + h_{d + 1 - i}$
%for  $0 \le i \le \left\lfloor \tfrac d 2 \right\rfloor $, %\]
%with at least one of these inequalities strict.
$h_{0} + \cdots + h_{i+1} \geq h_{d} + \cdots + h_{d - i}$ %\quad \textrm{ for } \quad 0 \le i \le \left\lfloor \tfrac d 2 \right\rfloor , \]
for $0 \le i \le \left\lfloor \tfrac d 2 \right\rfloor - 1$.
By Lemma \ref{Celtics} and the symmetry of the Eulerian numbers,
\begin{equation*} % \label{rhine}
h_{i + 1}(n) - h_{d - i}(n) = 2 \left( A(d-1,i + 1) - A(d - 1, i ) \right) g_{d-1} \, n^{d-1} + 2\sum_{r = 1}^{\lfloor \frac{d - 1}{2} \rfloor} \lambda_{d - 1 - 2r}(d) \, g_{d - 1 - 2r} \, n^{d - 1 - 2r},
\end{equation*}
for
$0 \le i \le   \lfloor \frac d 2 \rfloor - 1$, where $\lambda_{d - 1 - 2r}(d)$ is a function of $d$.
By Lemma \ref{deadset}, $g_{d - 1} \ge \frac{1}{2(d - 1)!}$ and hence by Theorem \ref{runner}, we can bound the ratios
$|g_{d - 1 - 2r}/g_{d  - 1}|$ in terms of $d$.  The above argument then shows that we can bound the absolute values of the roots of
$h_{i + 1}(n) - h_{d - i}(n)$ in terms of~$d$.
\end{proof}

%\begin{rem}
%Observe that the only part of the above theorem that does not follow from the nonnegativity of the coefficients of $h(t)$ is the fact that  $h_{i + 1}(n) > h_{d - i}(n)$ for
%$i = 0, \ldots, \left\lfloor \frac d 2 \right\rfloor - 1$ and $n \geq n_{d}$, which follows from the second assumption (by Remark \ref{Memphis}, this is the assumption that $a(t)$ is nonzero with  nonnegative  coefficients).
%\end{rem}

\begin{exa}
If $P$ is a $d$-dimensional lattice polytope and $h(t) = \delta_{P}(t)$, then $\U_{n} \delta_{P}(t) = \delta_{nP}(t)$ and the assumptions of the above theorem hold by Remark \ref{Memphis} and  Example \ref{JayHawk}. This establishes Corollary \ref{mainthm}.
\end{exa}

\section{Improving on the Bounds}\label{lowerboundsection}

One would like a bound on the integers $n_{d}$ and $m_{d}$ in Theorem \ref{Lebron}. In this direction, we will now show that for any positive integer $d$ and $n \geq d$, if $h(t)$ satisfies certain inequalities,  %the assumptions of Theorem \ref{ole},
then $h_{i + 1}(n) > h_{d - i}(n)$ for
$i = 0, \ldots, \left\lfloor \frac d 2 \right\rfloor - 1$.

We will continue with the notation of the previous section and consider a polynomial $h(t)$ %as in (\ref{star3}).
of degree at most $d$ with integer coefficients. By Lemma \ref{banana},
$h(t)$ has a unique decomposition
%\begin{equation*}
$h(t) = a(t) + b(t)$,
%\end{equation*}
where $a(t)$ and $b(t)$  are polynomials with integer coefficients satisfying
$a(t) = t^{d} \, a(\frac 1 t)$ and $b(t) = t^{d + 1} \, b(\frac 1 t)$.
Recall from Remark \ref{Memphis} that
the coefficients of $a(t)$ are strictly unimodal if and only if $h_{i + 1} > h_{d - i}$ for $0 \le i \le \lfloor \frac d 2 \rfloor - 1$.
By Lemma \ref{lions}, for any positive integer $n$,
%\[
$\U_{n}h(t) =  \E_{n} \left( h(t) \, (1 + t + \cdots + t^{n -1})^{d + 1} \right)$.
%\]
Setting
%\[
$\tilde{a}(t) := \E_{n} \left( a(t) \, (1 + t + \cdots + t^{n -1})^{d + 1} \right)$ and
%\]
%\[
$\tilde{b}(t) := \E_{n} \left(  b(t) \, (1 + t + \cdots + t^{n -1})^{d + 1} \right)$,
%\]
we have
%\[
$\U_{n}h(t) =  \tilde{a}(t) + \tilde{b}(t) \, $.
On the other hand, by Lemma \ref{banana}, we have a decomposition,
%we have a similar decomposition for $\delta_{nP}(t)$. We will write
%\begin{equation*}
$\U_{n}h(t) = a'(t) + \, b'(t) \, $,
%\end{equation*}
where $a'(t)$ and $b'(t)$  are polynomials with integer coefficients satisfying
$a'(t) = t^{d} \, a'(\frac 1 t)$ and $b'(t) = t^{d + 1} \, b'(\frac 1 t)$.
Our next goal is to express the polynomials $a'(t)$ and $b'(t)$ in terms of the polynomials
$\tilde{a}(t)$ and $\tilde{b}(t)$.
%$\U_{n} \left( a(t) \, (1 + t + \cdots + t^{n -1})^{d + 1} \right)$ and $\U_{n} \left( t \, b(t) \, (1 + t + \cdots + t^{n -1})^{d + 1} \right)$.
The next lemma says that $\tilde{b}(t)$ only contributes to $b'(t)$.

\begin{lemma}\label{Chicago}
%We can write $\E \left( tb(t) \, (1 + t)^{d + 1} \right) = t \, \tilde{b}(t)$ for some
The polynomial $\tilde{b}(t)$ satisfies
$\tilde{b}(t) = t^{d + 1} \, \tilde{b}(\frac 1 t)$.
\end{lemma}
\begin{proof}
If we let $f(t) = b(t)(1 + t + \cdots + t^{n -1})^{d + 1}$, then $f(t) = t^{n(d + 1)}f(\frac 1 t)$. Applying the operator $\E_{n}$ to both sides gives $\tilde{b}(t) = t^{d + 1}\tilde{b}(\frac 1 t)$.
\end{proof}

If we use the notation
\begin{equation}\label{Hally}
p(t) = a(t) (1 + t + \cdots + t^{n -1})^{d + 1} = \sum_{k = 0}^{n(d + 1) - 1} p_{k} \, t^{k},
\qquad \text{ then } \qquad
\tilde{a}(t) = \E_{n}p(t) \, .
%g_{k} = \sum_{i = 0}^{d} a_{i} \binom{d + 1}{k - i} \, .
\end{equation}
Observe that the symmetry of $a(t)$ implies that $p(t) = t^{n(d + 1) - 1}p(\frac 1 t)$. Hence we
may write
\[
p(t) = p_{0} + p_{1}t + \cdots + p_{n}t^{n} + \cdots + p_{n -1}t^{nd} + \cdots  + p_{1}t^{n(d + 1) - 2} + p_{0}t^{n(d + 1) - 1} ,
\]
which implies
\begin{equation}\label{UCLA}
\tilde{a}(t) = p_{0} + p_{n}t + p_{2n}t^{2} + \cdots + p_{\lfloor \frac d 2 \rfloor n}t^{\lfloor \frac d 2 \rfloor} +
p_{\lfloor \frac{ d + 1 }{ 2 } \rfloor n - 1} t^{\lfloor \frac d 2 \rfloor + 1} %+ g_{\lfloor \frac{ d + 1 }{ 2 } \rfloor n - 3} t^{\lfloor \frac d 2 \rfloor + 2}
+ \cdots + p_{2n -1}t^{d - 1} + p_{n-1}t^{d}.
\end{equation}
With the notation
$
a'(t) = a_{0}' + a_{1}'t + \cdots + a_{d}'t^{d},
$
we deduce the following lemma.

\begin{lemma}\label{Bird}
For $i = 0, \ldots, \left\lfloor \frac d 2 \right\rfloor$,
%\[
$a_{i}' = p_{0} + p_{n} + \cdots + p_{in} - p_{n - 1} - p_{2n - 1} - \cdots - p_{in - 1}$.
%\]
\end{lemma}
\begin{proof}
By Lemma \ref{Chicago}, to determine $a'(t)$ we only need to decompose $\tilde{a}(t)$ into its symmetric components as in Lemma \ref{banana}.
The result now follows from (\ref{coeff1}) and \eqref{UCLA}.
\end{proof}

If we fix $1 \leq k \leq \lfloor \frac{d}{2} \rfloor$, then Lemma \ref{Bird} implies that
%\[
$a_{k}' - a_{k - 1}' = p_{kn} - p_{kn - 1} \, $.
%\]
If $\gamma_{i}$ denotes the coefficient of $t^{i}$ in $(1 + t + \cdots + t^{n -1})^{d + 1}$, then,
%Since $g(t) = a(t) (1 + t + \cdots + t^{n -1})^{d + 1}$,
$p_{j} = \sum_{i = 0}^{d} a_{i}\gamma_{j - i} \,$ for $j = 0, \ldots, n(d + 1) - 1$,
and we conclude that
\begin{equation}\label{Xavier}
a_{k}' - a_{k - 1}' = \sum_{i = 0}^{d} a_{i}(\gamma_{kn - i} - \gamma_{kn - 1 - i}) \, .
\end{equation}

\begin{lemma}\label{Davidson}
The coefficients $\{\gamma_{i} \}$ %_{i = 0, \ldots, (n - 1)(d + 1) }$
of $(1 + t + \cdots + t^{n -1})^{d + 1}$
%If we write $(1 + t + \cdots + t^{n -1})^{d + 1} = \sum_{i = 0}^{(n - 1)(d + 1)} \gamma_{i} t^{i}$, then
are positive, symmetric and strictly unimodal.
\end{lemma}
\begin{proof}
%The result follows by induction on $d$.
%When $d = 1$,
%\[
%(1 + t + \cdots + t^{n -1})^{2} = 1 + 2t + \cdots + (n - 1)t^{n - 2} + nt^{n - 1} + (n - 1)t^{n} + \cdots + 2t^{2n -3} + t^{2n - 2},
%\]
%and the result holds.
The lemma follows from the fact that the product of two polynomials with positive, symmetric, unimodal coefficients has positive, symmetric, strictly unimodal coefficients.
%When $d > 1$, since $(1 + t + \cdots + t^{n -1})^{d + 1} =
\end{proof}

By Example \ref{JayHawk}, the assumption in the following lemma holds when $P$ is a $d$-dimensional lattice polytope and $h(t) = \delta_{P}(t)$.

\begin{lemma}\label{hoopdreams}
Suppose that $a(t)$ has positive integer coefficients and fix $1 \leq k \leq \lfloor \frac{d}{2} \rfloor$.
%For any $1 \leq k \leq \lfloor \frac{d}{2} \rfloor$, if
%If $n \geq \max{\left( \frac{d + 1}{d + 1 - 2k},\frac{d}{2k} \right)}$
If either $n$ and $d$ are even and $n \geq \frac{d}{d + 1 - 2k}$ or $n \geq \frac{d + 1}{d + 1 - 2k}$,
then $a_{k}' > a_{k - 1}'$.
\end{lemma}
\begin{proof}
%We need to show that the right hand side of (\ref{Xavier}) is positive.
%\[
%a_{k}' - a_{k - 1}' = \sum_{i = 0}^{d} a_{i}(\gamma_{kn - i} - \gamma_{kn - i - i}).
%\]
%By Theorem \ref{knock}, the coefficients $a_{i}$ are positive for $i = 0, \ldots, d$.
By Lemma \ref{Davidson},
if %$\left\lfloor \frac d 2 \right\rfloor \leq
$kn \leq \lceil \frac{(n - 1)(d + 1)}{2}  \rceil$ then $\gamma_{kn - i} - \gamma_{kn - 1 - i} \geq 0$ for $i = 0, \ldots, \left\lfloor \frac d 2 \right\rfloor$. Since the coefficients of $a(t)$ are positive, the right hand side of (\ref{Xavier}) is positive,
provided that $kn \leq \lceil \frac{(n - 1)(d + 1)}{2}  \rceil$. If $n$ and $d$ are even, the latter condition  holds if and only if
$n \geq \frac{d}{d + 1 - 2k}$. Otherwise, the condition holds if and only if $n \geq \frac{d + 1}{d + 1 - 2k}$.
%Since the coefficients of $a(t)$ are symmetric and nonnegative, it follows %from the strict unimodality of the $\gamma_{i}$
%that $a_{k}' > a_{k - 1}'$ when $d \le 2kn  \leq (n - 1)(d + 1)$. The latter condition holds if and only if $n \geq \max{\left( \frac{d + 1}{d + 1 - 2k},\frac{d}{2k} \right)}$.
\end{proof}

\begin{rem}
A similar lemma holds if we only assume that $a(t)$ is nonzero with nonnegative coefficients.
%If we only assume that $a(t)$ is nonzero with nonnegative coefficients, then we may manipulate the above proof to give a similar
%result to Lemma \ref{hoopdreams}.
%lower bound on $n$ to the one given in Lemma \ref{hoopdreams}.
\end{rem}

We now prove the main result of this section.

\begin{theorem}\label{Beijing}
Fix a positive integer $d$ and set $n_{d} = d$ if $d$ is even and $n_{d} = \frac{d + 1}{2}$ if $d$ is odd.  If $h(t)$ is a polynomial of degree  at most $d$
as in \eqref{star3} satisfying
%\[
$h_{0} + \cdots + h_{i+1} > h_{d} + \cdots + h_{d - i}$ %\quad \textrm{ for } \quad 0 \le i \le \left\lfloor \tfrac d 2 \right\rfloor , \]
for $0 \le i \le \left\lfloor \tfrac d 2 \right\rfloor - 1$,
%with at least one of the above inequalities strict,
then %for any $n \ge n_{d}$,
%If $P$ is a $d$-dimensional lattice polytope then, for any $n \geq d$,
%\[
$h_{i + 1}(n) > h_{d - i}(n)$ %\quad \textrm{ for } 0 \le i \le \left\lfloor \tfrac d 2 \right\rfloor - 1 \, .
%\]
for $0 \le i \le \left\lfloor \tfrac d 2 \right\rfloor - 1 \, $ and $n \ge n_{d}$.
\end{theorem}
\begin{proof}
By Remark \ref{Memphis}, we have assumed that $a(t)$ has positive integer coefficients and we need to show that the polynomial $a'(t)$ is strictly unimodal. The result now follows from
Lemma \ref{hoopdreams}.
\end{proof}

\begin{exa}\label{bigshot}
If $P$ is a $d$-dimensional lattice polytope and $h(t) = \delta_{P}(t)$, then the assumptions of the above theorem hold
by Example \ref{JayHawk}.
\end{exa}

%%%%%%%%%%%%%%%%%%%%%%%%%%%%%%%%%%%%%%%%%%%%%%%%%%%%%%%%%%%%%%%%%%%%%%%%%%%%%%%%%%

\section{Open Questions}\label{openquestionsection}

The main problem that remains concerns optimal choices (beyond Theorem \ref{Beijing}) for the integers $m_d$ and $n_d$ in Theorem \ref{Lebron} and Corollaries \ref{mainthm} and \ref{veronesecor}. We offer the following conjecture.

\begin{conj} 
If $\dim P = d$ then $\delta_{nP}(t)$ has distinct, negative real roots for $n \ge d$. 
\end{conj}

This conjecture holds for $d=2$, by the following argument:
A polynomial $1 + a_1 t + a_2 t^{2}$ has distinct real roots if and only if the discriminant 
$a_1^{2} - 4a_2 = (a_1 - 2)^{2} + 4(a_1 - a_2 - 1) > 0$. Hence the polynomial has distinct, real roots if 
$a_1 > a_2 + 1$.
If $\delta_{nP}(t) = 1 + h_1(n) \, t + h_2(n) \, t^{2}$, then $h_1(n) - h_2(n)$ equals the number of lattice points on the boundary of $nP$ minus 3 (see, e.g., \cite[Corollary 3.16 \& Exercise 4.7]{BRComputing}). If $n \ge 2$, each edge of $nP$ contains a lattice point that is not a vertex and hence $h_1(n) - h_2(n) \geq 3$. Thus $h_1(n) > h_2(n) + 1$ and $\delta_{nP}(t)$ has real roots. These roots have to be negative because the coefficients of $\delta_{nP}(t)$ are nonnegative.

Note that an example of a polytope $P$ such that $\delta_{P}(t)$ has complex roots is given by convex hull of $(0,1)$, $(1,0)$ and $(-1,-1)$, with $\delta_{P}(t) = 1 + t + t^{2}$.

%%%%%%%%%%%%%%%%%%%%%%%%%%%%%%%%%%%%%%%%%%%%%%%%%%%%%%%%%%%%%%%%%%%%%%%%%%%%%%%%%%

\bibliographystyle{amsplain}
%\bibliography{alan}
%\end{document}

\def\cprime{$'$}
\providecommand{\bysame}{\leavevmode\hbox to3em{\hrulefill}\thinspace}
\providecommand{\MR}{\relax\ifhmode\unskip\space\fi MR }
\providecommand{\MRhref}[2]{%
  \href{http://www.ams.org/mathscinet-getitem?mr=#1}{#2}
}
\providecommand{\href}[2]{#2}

\end{document}